\documentclass{amsart}
\usepackage[latin1]{inputenc}
\usepackage{amssymb,amsmath,latexsym}
\usepackage{color}
\usepackage[mathscr]{eucal}
\usepackage{enumerate}
\usepackage{amsthm}
\usepackage[all]{xy}

\def\begc{\begin{coro} }
\def\begd{\begin{defn} \em  }
\def\begeq{\begin{equation}}
\def\begeqas{\begin{eqnarray*}}
\def\begej{\begin{exmp} \em}
\def\begl{\begin{lem}}
\def\begp{\begin{prop} }
\def\begr{\begin{rem} \em }
\def\begt{\begin{thm} }

\def\calS{\mathscr{S}}
\def\cont{\subset}

\def\divUD{\frac{1}{2}}

\def\eleI{\ell^\infty}

\def\equS{:=}
\def\espV{\\[3pt]}

\def\infi{\infty}
\def\inn{\cap}

\def\NN{\mathbb{N}}

\def\nin{\notin}
\newcommand{\no}[1]{\|{#1}\|}

\def\Ome{\Omega}
\def\pt{\forall}
\newcommand{\raizc}[1]{\sqrt{{#1}} }
\def\RR{\mathbb{R}}

\def\setm{\setminus}
\def\Sig{\Sigma}
\def\sig{\sigma}
\def\span{\textrm{span}}
\def\subjun{_{j=1}^{n}}
\def\subkui{_{k=1}^{\infi}}
\def\subkuK{_{k=1}^{K}}
\def\submuM{_{m=1}^{M}}
\def\submuMU{_{m=1}^{M+1}}
\def\subnui{_{n=1}^{\infi}}

\def\subnuN{_{n=1}^{N}}

\def\terc{\end{coro}}
\def\terd{\end{defn}}
\def\terej{\end{exmp}}
\def\tereq{\end{equation}}
\def\tereqas{\end{eqnarray*}}
\def\terl{\end{lem}}
\def\terp{\end{prop} }
\def\terr{\end{rem} }
\def\tert{\end{thm}}
\def\tien{\rightarrow}

\def\un{\cup}
\def\vac{\phi}

\newtheorem{thm}{Theorem}[section]
\newtheorem{prop}[thm]{Proposition}
\newtheorem{lem}[thm]{Lemma}
\newtheorem{coro}[thm]{Corollary}
\newtheorem{rem}[thm]{Remark}
\newtheorem{defn}[thm]{Definition}
\newtheorem{exmp}[thm]{Example}

\numberwithin{equation}{section}

\date{\today}

\title[Equivalent Norms in a Banach Function Space]{Equivalent Norms in a Banach Function Space and the Subsequence Property}

\author[Calabuig]{J.M. Calabuig}
\address{Instituto Universitario de Matem\'{a}tica Pura y Aplicada
\\Universitat Polit\`{e}cnica de Val\`encia
\\
Camino de Vera s/n
\\
46022 Valencia
\\ Spain} \email{jmcalabu@mat.upv.es}

\author[Fern\'andez]{M. Fern\'andez Unzueta}
\address{Centro de Investigaci\'on en  Matem\'aticas, A.C.
\\
Jalisco S/N, Col. Valenciana CP: 36023 Guanajuato, Gto., M\'exico} \email{maite@cimat.mx}

\author[Galaz]{F. Galaz-Fontes}
\address{Centro de Investigaci\'on en Matem\'aticas, A.C.
\\
Jalisco S/N, Col. Valenciana CP: 36023 Guanajuato, Gto., M\'exico} \email{galaz@cimat.mx}

\author[S\'anchez]{E.A. S\'{a}nchez-P\'{e}rez}
\address{Instituto Universitario de Matem\'{a}tica Pura y Aplicada
\\ Universitat Polit\`{e}cnica de Val\`encia
\\
Camino de Vera s/n
\\ 46022 Valencia
\\ Spain} \email{easancpe@mat.upv.es}

\subjclass[2010]{46E30}

\keywords{Measure space, space of measurable functions, order, Banach function space.}

\thanks{All the authors were supported by Ministerio de Econom\'{\i}a, Industria y Competitividad (Spain). J.M. Calabuig and M. Fern\'andez  Unzueta under project  MTM2014-53009-P,
F. Galaz-Fontes under project MTM2009-14483-C02-01 and E.A. S\'anchez P\'erez  under project MTM2016-77054-C2-1-P}

\begin{document}

\begin{abstract}
Given a finite measure space $(\Ome,\Sig,\mu)$, we show that any Banach space $X(\mu)$
consisting of (equivalence classes of) real measurable functions defined on $\Ome$ such
that $f\chi_A  \in X(\mu) $ and $ \|f\chi_A \| \leq \|f\|,  \ \pt f \in X(\mu),
\ A \in \Sig$, and having the subsequence property, is in fact an ideal of measurable
functions and has an equivalent norm under which it is a Banach function space.  As an
application we characterize  norms  that are equivalent to a Banach function
space norm.
\end{abstract}

\maketitle

\section{Introduction and preliminaries}

 Throughout this paper all vector spaces we consider are real vector spaces and $(\Ome,
\Sig,\mu)$  is a finite measure space. If $X(\mu)$ is a normed (Banach) space
consisting of (equivalence classes of) real $\Sig$-measurable functions we will say that
$X(\mu)$ is a {\sl normed (Banach) space of measurable functions}.
A Banach space $X(\mu)$ is called  a  {\sl Banach function space} if it is an \textsl{ideal} of
$\Sig$-measurable functions and has a \textsl{Riesz norm}.  This means that if  $g$ is a
$\Sig$-measurable function,  $f \in X(\mu)$ and $| g| \leq |f|$,  then we always have
 $g \in X(\mu)$ and $\|g\| \leq \|f\|$.

In this paper (see Corollary \ref{Thm normaEFBEequiv}) we prove that a Banach space of
measurable functions $X(\mu)$ admits an equivalent norm $\|\cdot\|_V$ for  which
$(X(\mu), \|\cdot\|_V)$ is a Banach function space
if and only if it satisfies:

\begin{enumerate}
\item[(P1)]   There exists $C > 0$ such that for any $f \in X(\mu)$ we have
\[ f\chi_A  \in X(\mu), \ \pt A  \in \Sig, \textrm{ and } \|f\chi_A  \|_{X(\mu)} \leq C \|f
\|_{X(\mu)}.
\]

 \item[(P2)] 
 Any sequence that converges in $X(\mu)$ has a subsequence
converging pointwise $\mu$-a.e.
\end{enumerate}

In order to better understand the role of condition (P1) in the ideal property  of a
Banach function space, we first study the case for a Banach sequence space $\calS (X$) (see
definitions below), where $X$ is a Banach space with a Schauder basis. We show that
the unconditionality  of the basis corresponds, precisely, to condition (P1).

The vector space consisting of all real sequences $s = \{a_n\}$ is denoted
by $\calS$
and a {\it Banach sequence space} is simply a Banach space $X \cont \calS$. By taking $\Ome \equS
\NN, \Sig \equS 2^\NN$ and, for example, $\mu(A) \equS \sum _{n \in A} \frac{1}{n^2}, \ \pt A
\cont \NN$, we consider any Banach sequence space as a Banach space of measurable functions with
respect to a finite measure.

Given a Banach space $X$, recall that a sequence  $\{y_n\} \in X$ is {\it unconditionally
summable} if the series $\sum \subnui y_{\pi(n)}$ converges for each permutation
$\pi: \NN \tien \NN$.  A set $\{x_n\} \cont X$ is a {\it Schauder basis for $X$}, if for each
$x \in X$ there is a unique sequence $\{a_n \} \in \calS$ satisfying
\begeq \label{BaseS}
x = \sum \subnui a_n x_n.
\tereq

To each  Banach space $X$ with a Schauder basis $\{x_n\}$ we naturally associate the space
$\calS (X)$ formed by the sequences corresponding to the \lq coordinates' of its elements,
that is
\[
\calS (X) \equS \left\{  \{a_n\} \in \calS : \sum \subnui a_n x_n \in X \right\}.
\]
Moreover, we can define a linear map  $T: \calS(X) \tien X$ by
\[T(s) \equS \sum \subnui a_n x_n,\]
which is an isometric isomorphism   with the norm
$
 \|s\|_{\calS(X)} \equS \|T(s)\|_X.
$
In this way any Banach space $X$ with a Schauder basis can be
considered as a Banach sequence space, $\calS(X)$.
Additionally, by means of $T$  the order in $\calS(X) \cont \calS$  can be translated to
$X$. Namely, for $s, t \in \calS (X)$ we define
$
T(s) \leq T(t)  \textrm{ if }  s \leq t.
$

If $\{x_n\}$ is a Schauder basis for $X$ and, for each $x \in X$, the
sequence $\{a_n x_n\} $, where $\{a_n\}$ is as in (\ref{BaseS}), is unconditionally
summable, then $\{x_n\}$ is said to be an {\it unconditional basis for $X$}. In this case
the following properties are well known (\cite[p.~5, p.~344]{DJT}, \cite[pp.~15-16]{LT}):
For each   $\{a_n\} \in \calS (X)$
\begeq \label{MultLInfi}
  \sum \subnui c_n a_n x_n   \textrm{ converges for any } \{c_n\} \in \eleI,
\tereq
and there exists $C > 0$ such that
\begeq \label{MultLInfiCont}
\left\|\sum \subnui c_n a_n x_n \right\|_X \leq C \|t\|_\infi \left\|\sum \subnui a_n
x_n \right\|_X, \ \pt\, t = \{c_n\} \in \eleI.
\tereq
Of course, here  $\|\cdot\|_\infi$ is the norm of $\eleI$.

Before reformulating the above results, recall that for $s= \{a_n\} \in \calS$ and $A
\cont \NN$, then $s \chi_A \equS \{b_n\}$, where $b_n \equS a_n$ when $n \in A$ and $b_n
\equS 0$ when $n \nin  A$.

\begt \label{baseIncond}
Let $X$ be a Banach space. The following properties are equivalent:
\espV
{\rm i)} $\{x_n\}$ is an unconditional basis for $X$.
\espV
{\rm ii)} $\calS(X)$ is an ideal of $\calS$ and there exists $C > 0$ such that
\begeq \label{SXIdeal}
\|t\|_{\calS(X)} \leq C \|s\|_{\calS(X)}  \ \textrm{ if } \  t, s \in \calS(X)  \textrm{ and }  |t|
\leq |s|.
\tereq

\noindent {\rm iii)}  There exists $C > 0$ such that for any $s \in \calS(X)$ we have
{\em \begeq \label{NormaSX}
s\chi_A  \in \calS (X), \ \pt A \cont \NN, \mbox{ and } \|s\chi_A  \|_{ \calS(X)} \leq C \|s\|_{ \calS(X)}.
\tereq }
\tert
\proof  Since the other implications are well known, we will only prove that i) implies ii).
Take $s =\{a_n\} ,  t= \{b_n\} \in \calS$  such that $ |t| \leq |s|$. Let us define
\[
c_n=\begin{cases}
\displaystyle \frac{b_n}{a_n}& \ a_n\neq 0\\
0 &a_n=0
\end{cases}\]
and consider  $r \equS \{c_n\}$. Assume that $s \in \calS(X)$, so that
 the sequence $\{a_n x_n\}$ is unconditionally summable. Since $r \in \eleI$, applying
(\ref{MultLInfi})  we find that $t = rs \in \calS(X)$. This shows that $\calS (X)$ is
an ideal of $\calS$.

From (\ref{MultLInfiCont}) we also have $\big\|T(r  s)\big\|_X \leq  C \|r\|_{\infi} \big\|T (s)\big\|_X$.
This corresponds to
\[
 \|t\|_{\calS(X)} =  \big\| r s \big\|_{\calS (X)} \leq C \|r\|_\infi \|s \|_{\calS(X)} \leq C \|s\|_{\calS (X)}.
\]
\qed

\

Let us consider  now  the general case of  a Banach  space of measurable functions, say $X(\mu)$.
Although  condition i) in Theorem \ref{baseIncond} may have no  sense when we consider   $\Sig$
instead of $2^\NN$ and  $X(\mu)$ instead of $\calS(X)$, conditions ii) and iii)  still  have sense,
as  the example of the  sequence space $\eleI$  shows: it  has properties ii) and iii) but, since
it is not separable, it does not have a Schauder basis.

Then,  it is a  natural question to ask if  conditions ii) and iii) remain
equivalent in this more general context.

In this paper we  address  this  question.  Actually, we prove that ii) is equivalent to
 condition (P1) -which would correspond to condition iii) in Theorem \ref{baseIncond}- plus the
 already introduced subsequence property (P2).

The paper is organized as follows:
Section \ref{sec: prelim}  collects the preliminary results we will need further, while
Section \ref{sec: subseq} is devoted to study what we call {\sl subsequence property}, namely,
condition (P2) in the Introduction. Condition (P1) is studied in Section \ref{sec: rectangular}.
We named \emph{rectangular} those spaces $X(\mu)$, or norms, that satisfy (P1). It is of interest
to note that this property does not require explicitly that $X(\mu)$ be an ideal (see Section
\ref{sec: rectangular}).

In Section \ref{sec: main} we prove the main result of the paper, namely Theorem \ref{equivEFB}.

As a consequence, as we note in Section \ref{sec: conseq}, it follows readily that   basic
properties of Banach functions spaces  remain valid in the more general context of Banach rectangular spaces
with the subsequence property.   In Section \ref{sec: conseq} we also show that
rectangular sequence spaces always have the subsequence property (see Proposition
\ref{propRectSeq}), that is, in this case (P1) implies (P2). This result explains why in the case
of sequence spaces, condition ii) in Theorem \ref{baseIncond}
is equivalent to  iii) (namely to (P1)).  Observe, also, that
our  setting includes those Banach spaces $X$ with unconditional basis.

\section{Preliminaries}\label{sec: prelim}

Let $(\Ome,\Sig,\mu)$ be a finite measure space. Then
$L^0(\Sig)$ consists of all $\Sig$-measurable functions $f: \Ome \tien \RR$, $S(\Sig)$ is
the vector space formed by the simple functions $s: \Ome \tien \RR$ and
$L^0(\mu)$ is the vector space consisting of the equivalence classes of $\Sig$-measurable
functions, where two functions are equivalent when $f = g\, \mu$-a.e. As usual,
$\|\cdot\|_\infi$  denotes  the norm in $L^\infi(\mu)$.

In $L^0(\mu)$ we consider the canonical order, given by $f \leq g$ if $f \leq g\,
\mu$-a.e. Given $A \cont L^0 (\mu)$, then $A^+ \equS \{f \in A: f \geq 0 \}$. The notation
$X(\mu)$ will indicate that  $X(\mu)$ is a vector subspace of $L^0 (\mu)$; we define
$S(X(\mu)) \equS S(\Sig) \inn X(\mu)$ and $X(\mu)^+ $ is the \emph{positive cone} of
$X(\mu)$.

A space $X(\mu)$ is a \emph{Riesz space} if $\max\{f,g\} \in X(\mu), \ \pt f, g
\in X(\mu)$ \cite{ZaIOT}.  Since
\begeq \label{maximo}
\max\{f,g\} = \frac{f - g + |f - g|}{2}, \ \pt f, g \in L^0(\mu),
\tereq
if the absolute value $|h| \in X(\mu)$ for any $h \in X(\mu)$, then $X(\mu)$ is a Riesz
space.

If $f \in X(\mu)$  when  $f \in L^0(\mu), g \in X(\mu)$  and $0 \leq  |f | \leq
|g|$, the  vector space $X(\mu)$ is an \emph{ideal} of $L^0(\mu)$. Notice that in this
case, given $f \in L^0 (\mu)$, we have $f \in X(\mu)$ if and only  if  $|f| \in X(\mu)$.
Hence an ideal is always a Riesz space.

A norm $\|\cdot\|$ on a Riesz space $X(\mu)$ is a \emph{Riesz norm}
\cite[p.~48, p.~438]{LuZa}, (or $X(\mu)$ is a \emph{normed Riesz space}) when $f, g
\in X(\mu)$ and $|f| \leq |g|$ imply that $\|\,f\,\| \leq \|\, g\,\|$. In this case we
clearly have $\no{f} = \|\,|f|\,\|, \ \pt f \in X(\mu)$.

A normed Riesz space $X(\mu)$ is a \emph{$\mu$-normed function space} ($\mu$-n.f.s. or n.f.s.
for short), when $X(\mu)$ is an ideal in $L^0(\mu)$. If additionally $X(\mu)$ is  complete,
then it is a Banach function space (Banach f.s.).

\subsection{The subsequence property}\label{sec: subseq}

In this subsection we study the subsequence property and relate it to other basic concepts
connected with the order. Although the subsequence property is a well-known feature of Banach function spaces, there are some current
research papers that study natural generalizations for abstract Banach lattices  (see for example \cite[Proposition 4.1]{(a)}, Corollary 3.5 and the references in this paper).
However, in our case we work in the setting of spaces of measurable functions, and so we can fix an underlying structure of measurability, and the order is naturally  defined as the a.e. order.

Recall that a  normed space $X(\mu)$ has the \emph{subsequence property} if for any
sequence  $\{f_n\} \cont X(\mu)$ and $ f \in X(\mu)$ such that  $f_n \tien f$, there is a
subsequence $\{f_{n(k)}\}$ converging pointwise to $f\, \mu$-a.e.

It is well known that the subsequence property can be characterized as follows \cite[18.14]{ZaI}

\begp \label{prop continuous inclusion}
A  normed space $X(\mu)$ has the \emph{subsequence property} if and only if the inclusion
map $X(\mu) \subset L^0 (\mu)$ is continuous.
\terp

Given a Riesz space $X(\mu)$,  then $X(\mu)^+$  is its  positive cone and  $V$
is  the  function
\[
V: X(\mu) \tien X(\mu), \textrm{ defined by } V(f) \equS |f|.
\]

The next lemma relates  the above concepts on Riesz spaces with a norm.

\begl  \label{VContinuaRiesz} Let $X(\mu)$ be a Riesz space with a norm.
\espV
i) If  $V$ is continuous, then  $X(\mu)^+$ is closed.
 \espV
ii) If $X(\mu)$ has the subsequence property, then $X(\mu)^+$ is closed.
\terl

\noindent
Proof.
 i) Let $\{f_n\} \cont X(\mu)^+$ and $f \in X(\mu)$ be such that $f_n \tien f$. Since
$V$ is continuous, this implies $f_n = |f_n| \tien |f|$. Hence $f = |f| \in X(\mu)^+$.

 ii) Let $\{f_n\} \cont X(\mu)^+$ and $f \in X(\mu)$ be such that $f_n \tien f$.
Next, take  a subsequence $\{f_{n(k)}\}$ such that $\{f_{n(k)}\}$ converges to $f \,
\mu$-a.e. Since $f_{n(k)} \geq 0 \,\mu$-a.e., it follows that $f \geq 0\, \mu$-a.e., that
is $f \in X(\mu)^+$.
\qed \vspace{6pt}

The following   normed sequence space  is an example of  a Riesz space whose positive cone
is not closed. Consequently, it  does not have the subsequence property and  the function
$V$ is not continuous on it.

\begej For each $n \in \NN$, let $v_n \equS e_1 + \cdots + e_n$, where $\{e_n\}$ is the
canonical orthonormal basis in $(\ell^2,\|\cdot\|_2)$. Then we define
\[
X(\mu) \equS \span\{v_n:n \in \NN \}, \ \ \left\|\sum \subnuN a_n v_n \right\| \equS
\left\|\sum \subnuN \frac{ a_n}{n} e_n \right\|_2.
\]
 Then $\|\cdot\|$ is a norm in $X(\mu)$ and $\|v_n\| = \frac{1}{n} \tien 0$.

Recall that we consider in $X(\mu)$ the natural order between sequences.  Then  $-e_1 + v_n \in
X(\mu)^+, \ \pt n \in \NN$, and $-e_1 + v_n \tien -e_1 \nin X(\mu)^+$. This shows that
$X(\mu)^+$ is not closed.

Since $|f| \in X(\mu), \ \pt f \in X(\mu)$, and so $X(\mu)$ is a Riesz space,
it follows from  i) of Lemma \ref{VContinuaRiesz} that the function $V$ is not  continuous.
\terej

Even when   $X$ is both a Riesz space and a Banach sequence space with the subsequence
property,  the  function  $V$ defined by $V(s) \equS |s|$ can be  discontinuous, as an
example of D.  H. Fremlin shows \cite[2XD]{Fr}.

Let us finish this section by recalling  some relations of the subsequence property with other classical and recently introduced Banach lattice properties, specially those lattices that
are endowed with a natural structure of measurability, as the Banach function spaces. In this paper, a particular type of monotonicity for the norm of the function spaces involved will be relevant. Monotonicity is often related to order continuity or Fatou type properties in Banach function spaces, and plays a relevant role regarding the relation among convergence in measure and norm convergence providing  a sort of converse of the subsequence property: every increasing order bounded sequence converges in the norm.
 Thus, several concepts associated to monotonicity properties for Banach lattices are defined in terms of sequences and are
naturally connected with the convergence in measure ---as the subsequence property--- in the case of Banach function spaces (see for example Theorem 2.5 of \cite{ForaHud}). More results for concrete Banach function spaces can be found in the excellent survey  \cite{ForaHud}. For current developments on the role of sequential properties and the geometry of Banach function spaces see \cite{Ciel}, and \cite{Hud} for the vector valued case.

\subsection{Rectangular Function Spaces}\label{sec: rectangular}
We  present in what follows the main definition of this paper and the fundamental properties of the spaces satisfying it. As the reader may notice, our definition lies in a basic monotonicity property as the ones commented in the previous section for a certain class of functions: an inequality in the almost everywhere order implies an inequality for the corresponding norms.

\begd A normed space of measurable functions $X(\mu)$ is a \emph{ rectangular function
space} (rectangular f.s.) if there exists $C > 0$ such that
\begeq \label{dRectangular}
\chi_A f \in X(\mu) \textrm{ and } \|\chi_A f \| \leq C \|f\|,  \ \ \pt f \in X(\mu), A \in
\Sig.
\tereq
In this situation, the norm will also be called \emph{rectangular}.
When $C = 1$ we will say that $X(\mu)$, and its norm, is \emph{strictly rectangular}.
If a rectangular (or strictly rectangular) f.s. is also norm complete, we will say it is a
\emph{Banach rectangular (or strictly rectangular) f.s}.
\terd

\begej \label{prop BfsRectang} Every normed function space $X(\mu)$  is a strictly
rectangular f.s. and if $\|\cdot\|_0$ is an equivalent norm for $X(\mu)$, then
$\|\cdot\|_0$ is rectangular.
\terej
\proof
Given a n.f.s $X(\mu)$, consider $f \in X(\mu), \ A \in \Sig$. Since $|\chi_A f|
\leq |f|$, we have $\chi_A f \in X(\mu)$ and
\[
\| \chi_A f \| = \big\| \chi_A |f|\,\big\| \leq \big\|\,|f|\,\big\| = \| f\|.
\]
This shows that  $X(\mu)$ is a strictly rectangular f.s.

Next, let us take $0 < a < b$ such that $a \|f\|_0 \leq \|f\| \leq b \|f\|_0$.
Then, for $A \in \Sig$ we have $a \|\chi_A f \|_0 \leq  \|\chi_A f\| \leq  \|f\|
\leq  b \|f\|_0$. It follows that the norm $\|\cdot\|_0$ is rectangular.
\qed

\

If $X(\mu)$ is a rectangular f.s., we next show that $|f| \in X(\mu)$ when $f \in X(\mu)$.
By (\ref{maximo}) it follows that any rectangular f.s. is a Riesz vector space. \vspace{6pt}

 Given $f \in L^0 (\mu)$, take $A \equS \{ x \in \Ome: f(x) \geq 0 \}$ and $B \equS \Ome
 \setm A$. Then for the positive and negative parts of $f$ we have
\begeq \label{fMas}
 f^+ = f \chi_A = |f| \chi_A, \ \ f^- = - f \chi_B =  |f| \chi_B.
\tereq
Since $f = f^+ - f^-$ and $|f| = f^+ + f^-$, we obtain immediately the following.

\begl  \label{equivalentes} Let  $X(\mu)$ be a rectangular f.s.
\espV
i)  Let $f \in L^0(\mu)$.  Then the following properties are equivalent:

\ \ \ a) $f \in X(\mu)$.

\ \ \ b) $f^+, \ f^- \in X(\mu)$.

\ \ \ c) $|f| \in X(\mu)$.
\espV
ii) $\divUD \|f\| \leq \big\|\,|f|\,\big\| \leq 2 \|f\|, \ \pt f \in X(\mu)$.
\terl

 Our next example  shows  that even in  $\RR^2$ there are  rectangular
 norms that are not Riesz norms.

\begej Let
\[
B \equS \{(x,y) \in \RR^2: xy \geq 0, x^2 + y^2 \leq 1 \} \un \{(x,y) \in \RR^2:
x y \leq 0,  |y - x| \leq 1\}.
\]
Then $B$ is an absorbing closed convex set  satisfying $-B = B$. Hence, by means of its
Minkowski functional we obtain a norm $\|\cdot\|_B$  for $\RR^2$, whose closed unit ball
is $B$ (\cite{Ru}, p. 25). The norm $\|\cdot \|_B$ is
rectangular. Take now $f \equS (-1,1)$. Then  $\|\,|f|\, |_B = \|(1,1)\| = \raizc{2}$.
On the other side we have $\|f\| = 2$. Hence $\no{f} > \big\|\,|f|\,\big\|$ and so this norm is
not a  Riesz norm.
\terej

 The next result proves that  for  Banach  rectangular f.s., implication  ii) in  Lemma
 \ref{VContinuaRiesz} is, indeed,   an equivalence.

\begl  \label{lEquiv recta bis} Let   $X(\mu)$ be a Banach rectangular f.s. Then, the
following properties are equivalent:
\espV
i) $X(\mu)$ has the subsequence property.
\\
ii) $ X(\mu) \subset L^0(\mu)$ continuously.
\\
iii)   $X(\mu)^+$ is closed.
\terl

\noindent
Proof. The equivalence between i) and ii) is already proved in \ref{prop continuous inclusion}.
That i) implies iii) is proved  in Lemma \ref{VContinuaRiesz}. To prove the converse
implication, consider $\{f_n\} \cont X(\mu)$ and $f \in X(\mu)$ such that $f_n \tien
 f$. Given $k \in \NN$ let us take $n(k) \in \NN$ so that $\|f - f_{n(k)}\| \leq 2^{-k}$
and $n(k) < n(k+1)$. Since $X(\mu)$ is a Banach space, using Lemma \ref{equivalentes} we
conclude that $S \equS \sum \subkui |f - f_{n(k)}|  \in X(\mu)$. Define $s_K \equS \sum
\subkuK |f - f_{n(k)}|,  K \in \NN$. For a fixed $K \in \NN$ and $\ell \in \NN$,
we have $s_{K + \ell} \geq  s_K$.  Letting now $\ell \tien \infi$ and taking into account
that $X(\mu)^+$ is closed, it follows that $S \geq s_K, \ \pt K \in \NN$. From here we can
find a set $A \in \Sig$ such that $\mu(A) = 0$ and $S(x) \geq S_K(x), \ \pt x \in \Ome
\setm A, \ K \in \NN$. This implies that the series $\sum \subkui |f - f_{n(k)}|$ converges
in $A$ and so $f_{n(k)} \tien f$ in $A$.
\qed \vspace{6pt}

It will follow from Theorem 5.8 that for a Banach rectangular f.s. the continuity  of $V$
is also equivalent with the properties considered in the above lemma.

\begd \label{Monotona}
The norm in $X(\mu)$ is \emph{monotone} if there exists $C > 0$
such that for $f, g \in X(\mu)^+$ satisfying $f \leq g$, we have $\no{f} \leq   C \no{g}$.
When $C =1$, we will say the norm is \emph{strictly monotone}.
\terd

\begl  \label{VContinua}
If  $X(\mu)$ is a rectangular f.s. with a monotone norm, then  the function $V$ is
continuous. More precisely,
\begeq \label{VLip}
\big\| \,|g| - |h|\, \big\| \leq 4 C \|g - h \|, \pt g, h \in X(\mu).
\tereq
\terl
\noindent
Proof. Take $C > 0$ as in Definition \ref{Monotona}. Let $g, h \in X(\mu)$. Then
\[
0 \leq \big|\, |g| - |h|\, \big| \leq | g -h|.
\]
Using now Lemma \ref{equivalentes} and the hypothesis we obtain
\[
\divUD \big\| \, |g| - |h|\, \big\| \leq \left\|\, \big||g| - |h|\big|\,\right\| \leq C \big\| \,|g - h|\,\big\| \leq 2C
  \| g - h \|,
\]
and the conclusion follows. \qed

Using Lemmas \ref{VContinua}, \ref{VContinuaRiesz} and \ref{lEquiv recta bis}
we obtain the following well known result.

\begc \label{BfsSubseqPpty}
Any Banach function space $X(\mu)$ has the subsequence property.
\terc

\section{Associated Banach Function Space norm}\label{sec: main}

In this section we prove the main result of this paper, namely Theorem \ref{equivEFB}, which
gives a   characterization of Banach spaces of measurable functions admitting an equivalent
norm  with which they are  Banach function spaces.

\begd Let $X(\mu)$ be a rectangular f.s. Given a function $f \in L^0(\mu)$, let $\Sig_f \equS
\{A \in \Sig:  f \chi_A \in X(\mu)\}$ and $S(f,X(\mu)) \equS \span \{ f\chi_A: A \in \Sig_f\}
\cont X(\mu)$.
\terd

Note that $\Sig_f$ is always a ring of subsets, that is, $\vac \in \Sig_f$ and $A \un B,
A \setm B \in \Sig_f$, when $A, B \in \Sig_f$. We also have $S(\chi_\Ome,X(\mu)) = S(X(\mu))$
and $\Sig_f = \Sig, \ \pt f \in X(\mu)$.

\begl Let $X(\mu)$ be a rectangular f.s. If $g \in S(f,X(\mu))$, then it can be expressed as
$g = \sum \subnuN a_n f \chi_{A_n} $, where $\{A_n\}_{n= 1}^{N}\cont \Sig_f$ consists of
disjoint sets.
\terl

\noindent
Proof. Let $g \in S(f,X(\mu))$. Then $g = \sum \submuM b_m f \chi_{B_m}  $ where $\{B_m\}_{m= 1}
^{M} \cont \Sig_f$. We will proceed by induction on $M$. For $M = 1$, the conclusion is clear.

Let us now consider $g = \sum \submuMU b_m f \chi_{B_m}$ where $\{B_m\}_{m= 1}^{M+1} \cont
\Sig_f$. By the induction hypothesis there is a finite set $J$ and a collection of disjoint
sets $\{C_j: j \in J \} \cont \Sig_f$ such that $\sum \submuM b_m f \chi_{B_m} = \sum _{j \in J}
c_j f \chi_{C_j}$. If $B_{M +1} \inn  C_j = \vac, \pt j \in J$, or $B_{M+1} = C_j$ for some $j
\in J$, the conclusion is clear. Assume this is not so, take   $J_0 \equS \{j \in J:
B_{M +1} \inn C_j \neq \vac \}$ and define $D_{M + 1} \equS B_{M+1} \setm \un _{j \in J}  C_j
$.
Then
\begeq
D_{M +1} \in \Sig_f, \ \ D_{M + 1} \inn C_j = \vac, \ \pt j \in J.
\tereq

Now  for each $j \in J_0$, let $D_j \equS C_j \inn B_{M+1}$ and $E_j \equS  C_j \setm
B_{M+1}$ Then $D_j, E_j \in \Sig_f, \ D_j \inn E_j = \vac$ and $C_j = D_j  \un E_j$.  We
also have
\begeq
B_{M+1} \inn (\un _{j \in J}  C_j)  =  \un _{j \in J_0} D_j.
\tereq

Hence
\begeqas
\sum \submuMU b_m f \chi_{B_m} & = & \sum _{j \in J} c_j f \chi_{C_j} + b_{M+1} f
\chi_{b_{M+1}}
\\
 & = & \sum _{j \in J \setm J_0} c_j f \chi_{C_j} + \sum _{j \in J_0} c_j f \chi_{C_j}
+ b_{M+1}
\left(\sum_{j \in J_0} f \chi_{D_j} + f\chi_{D_{M+1}} \right)
\\
\phantom{XXX} =  && \hspace{-24pt}  \sum _{j \in J \setm J_0}  c_j  f \chi_{C_j} +
\sum _{j \in  J_0} c_j
 f \chi_{E_j} + \sum _{j \in J_0} (c_j  + b_{M+1}) f \chi_{D_j} + b_{M+1} f
 \chi_{D_{M+1}},
\tereqas
and the conclusion follows. \qed

\begl Let $X(\mu)$ be a strictly rectangular f.s.  Take $f \in L^0(\mu)^+$ and  a finite collection
of disjoint  sets $\{A_j\}_{j= 1}^{n} \cont \Sig_f$. If  $0 \leq a_j \leq b_j,j = 1, \ldots,
n,$ then
\[
\left \|\sum \subjun a_j \chi_{A_j}f \right\| \leq \left\|\sum \subjun b_j \chi_{A_j} f \right\|.
\]
\terl
\proof  Let $n \in \NN$ be given. For $n = 1$ the conclusion is clear. So we now consider
 $n \geq 2$ and first assume that $a_1 = b_1, \ldots, a_{n-1} = b_{n-1}$ and $0 \leq a_n
< b_n$.   Take $ x \equS \sum _{1 \leq j < n} a_j \chi_{A_j} f$ and
  $A \equS  \un_{1 \leq j < n} A_j$. Since the norm is rectangular, we have
 \[
  \| x\| =  \|(x + b_n \chi_{A_n} f)  \chi_ A\| \leq  \|x + b_n \chi_{A_n}f \|.
  \]
Then
 \begeqas
\big\|(x + a_n \chi_{A_n}) f \big\| & = & \left \|\frac{a_n}{b_n}\left(  x + b_n \chi_{A_n}f \right)
+ (1 - \frac{a_n}{b_n}) x \right \|
\\[3pt]
& \leq & \frac{a_n}{b_n}\, \| x + b_n \chi_{A_n}f \| + \left(1 - \frac{a_n}{b_n} \right) \|x\|
\\[3pt]
& \leq  &\| x + b_n \chi_{A_n} f\|.
\tereqas

The general case is now obtained by applying what we have just proved to $ \sum _{1 \leq j
< n} a_j \chi_{A_j} f + a_n \chi_{A_n}f$ and $ \sum _{1 \leq j < n} a_j \chi_{A_j}f + b_n
\chi_{A_n}f$, then  to $ \sum _{1 \leq j < n-1} a_j \chi_{A_j}f$ $+ a_{n-1} \chi_{A_{n-1}}f +
b_n \chi_{A_n}f$ and $ \sum _{1 \leq j < n-1} a_j \chi_{A_j}f + b_{n-1}
\chi_{A_{n-1}}f + b_n \chi_{A_n}f$, until we finish with $ a_ 1\chi_{A_1}f +   \sum _{1 < j
\leq  n} b_j \chi_{A_j}f$ and  $ b_ 1\chi_{A_1}f +   \sum _{1 < j \leq  n} b_j \chi_{A_j} f$.
\qed \vspace{6pt}

\begp \label{nMonotona}
The norm in a strictly rectangular f.s. $X(\mu)$ is always strictly monotone on $S(f,X(\mu)),
\ \pt f \in L^0(\mu)$, and so
\begeq \label{VLipD}
\big\| \,|g| - |h|\, \big\| \leq  4 \|g - h\|,   \ \pt g, h \in S(f,X(\mu)),
\tereq

If $f \in X(\mu)$ we additionally have
\begeq \label{producto}
\big\|\,|sg|\, \big\|  \leq  \|s\|_\infi \big\|\,|g|\,\big\|,  \ \ \|s g \|  \leq  4 \|s\|_\infi
\|g\|,  \ \ \pt s \in S(\Sig), \ g \in S(f,X(\mu)).
\tereq
\terp

\noindent
Proof. Consider  $g, h \in S(f,X(\mu))$ such that $0 \leq g \leq h$ and let us express
$g = \sum \subjun a_j \chi_{A_j}f$ and $s = \sum \subjun b_k \chi
_{B_j}f$, where $\{A_j\}, \{B_k\} \cont \Sig $ are finite partitions of $\Ome$. Then
$g = \sum_{j,k} a_j \chi_{A_j \inn B_k} f$ and $h = \sum_{j,k} b_k \chi_{A_j \inn B_k}f$.
If $A_j \inn B_k$ has positive measure then $a_j \leq b_k$ and if $A_j \inn B_k$ has
measure zero we can go on by eliminating this term. Proceeding in this way we can apply the
above lemma to obtain that $\|g\| \leq \|h\|$.

Inequality (\ref{VLipD}) follows now from (\ref{VLip}) in Lemma \ref{VContinua}.

Assume now that $f \in X(\mu)$. Let $s \in S(\Sig)$ and $g \in S(f,X(\mu))$.  Since
$\chi_A f \in X(\mu), \pt A \in \Sig$, it follows that $s g \in X$ and $|s g | \leq
\|s\|_\infi |g|$. Since the norm is monotone on $S(f,X(\mu))$, we obtain the first
inequality in (\ref{producto}). Using now Lemma \ref{equivalentes} we have
\[
\|s g \| \leq 2 \big\|\,|s g|\,\big\| \leq 2 \|s\|_\infi \big\|\, |g|\,\big\| \leq 4
\|s\|_\infi \|g\|. \qed \vspace{6pt}
\]

\

In the sequel, we will consider the spaces to be complete. The next example shows that
completeness is basic for $X(\mu)$ to be an ideal.

\begej
Let $\Ome \cont \RR^n$ be a (Lebesgue) measurable set having positive measure and take
$\mu$ as Lebesgue's measure defined on the $\sig$-algebra of measurable subsets of $\Ome$.
Consider   $X(\mu) \equS S(\Sig)$  with the corresponding  $L^1$ norm. Then $X(\mu)$ is
a rectangular f.s. and is not an ideal in $L^0 (\mu)$.
\terej

\begl \label{lSubseqProp} Let $X(\mu)$ be a Banach space of measurable functions with the
subsequence property. If $f \in L^0 (\mu)$ and $\{f_n\} \cont X(\mu)$ is a Cauchy sequence
that has a subsequence $\{f_{n(k)}\}$ which converges pointwise $\mu$-a.e. to $f$, then
$f \in X(\mu)$, and $f_n \tien f$ in $X(\mu)$.
\terl

\noindent
Proof.  Let $g_k \equS f_{n(k)}, \ \pt k \in \NN$. Then $\{g_k\} \cont X(\mu)$  is also a
Cauchy sequence and so  there is some $g \in X(\mu)$ such that $g_k \tien g$. To obtain
the conclusion we will show that $g = f\,\mu$-a.e. Since $X(\mu)$ has the
subsequence property, let $\{g_{k(j)}\}_j$ be a subsequence such that $g_{k(j)} \tien g$
when $j \tien \infi$, pointwise $\mu$-a.e. By the hypothesis we also have $g_{k(j)} \tien
 f$ pointwise $\mu$-a.e. it follows that $f = g\,\mu$-a.e. \qed \vspace{6pt}

\begt \label{equivEFB} Let  $X(\mu)$ be a Banach strictly rectangular f.s. with the subsequence
property.
\espV
i) If $h \in L^\infi(\mu)$ and $f \in X(\mu)$, then $hf  \in X(\mu)$ and
\[
\ \ \ \big\|\,|hf|\,\big\| \leq \big\|\,|h|\,\big\|_\infi \big\|\,|f|\,\big\| , \ \ \  \|h f\| \leq 4 \|h\|_\infi
\|f\|.
\]
ii) $X(\mu)$ is an ideal in $L^0(\mu)$.
\espV
iii) The norm in $X(\mu)$ is monotone, (\ref{VLip}) holds and the function $V$ is continuous.
\espV
iv) The function
\begeq \label{normaV}
\|f\|_V \equS \big\|\,|f|\,\big\|, \ \pt f \in X(\mu).
\tereq
defines an equivalent norm for $X(\mu)$. Moreover, the norm $\|\cdot\|_V$ coincides with that
of $X(\mu)$ in $X(\mu)^+$ and, with the norm $\|\cdot\|_V$, the space $X(\mu)$ is a Banach f.s.
\tert

\noindent
Proof. i) Take $h \in L^\infi(\mu), \ f \in X(\mu)$ and let $\{s_n\} \cont S(\Sig)$,
be such  that $s_n \tien h$ in $L^\infi(\mu)$. Since $f \in S(f,X(\mu))$, it now follows
from the second inequality in (\ref{producto}) that $\{s_n f\} \cont X(\mu)$ is a Cauchy
sequence.

On the other side,  we can find a subsequence $\{s_{n(k)}\}$  such that $s_{n(k)} \tien h$
pointwise $\mu$-a.e. Then  $s_{n(k)} f \tien hf$ pointwise $\mu$-a.e. This allows us to
apply Lemma \ref{lSubseqProp} and conclude that $hf \in X(\mu)$ and $ s_n f \tien hf $ in
$X(\mu)$. Hence from (\ref{producto}) we obtain
\[
\|hf\| = \lim \|s_n f\| \leq 4 \lim \|s_n\|_\infi  \|f\| = 4 \|h\|_\infi \|f\|.
\]

Since  $s_n \tien h$ in $L^\infi(\mu)$, it follows that $|s_n| \tien |h|$ in $L^\infi(\mu)$.
As above, it follows that this implies $ |s_n| |f| \tien |h||f| $ in $X(\mu)$. Using now the
first inequality in (\ref{producto}) we have $\big\|\,|hf|\,\big\| \leq \big\|\,|h|\,\big\|_\infi \big\|\,|f|\,\big\|$.

ii) Take $f \in L^0 (\mu)$ and $g \in X(\mu)$ such that $|f| \leq |g|$. By Lemma
\ref{equivalentes} it is enough to show that $|f| \in X(\mu)$. For this, note that we can express
$|f| = h |g|$ for some $h \in L^\infi(\mu)$ with $\|h\|_\infi \leq 1$. Since $|g| \in X(\mu)$,
applying i) we obtain that $|f| \in X(\mu)$.

iii) Let $f, g \in X(\mu)$ be such that $0 \leq f \leq g$. Then we can express $f = h g$ where
$h \in L^\infi(\mu)$ and $\|h\|_{L^\infi(\mu)} \leq 1$. By i) we now have $\|f\| = \|hg\| \leq
\|g\|$.

iv) The other properties of a norm being clear, we will only establish the triangle inequality.
So, take $f, g \in X(\mu)$. Then $0 \leq |f + g| \leq |f| + |g|$. Since the norm in $X(\mu)$
is monotone we have
\[
\big\|\, |f + g|\, \big\| \leq \big\|\, |f| + |g|\, \big\| \leq \big\|\,|f|\,\big\|  + \big\|\,|g|\,\big\|.
\]

This shows that $ \|f +g \|_V \leq \|f\|_V + \| g\|_V, \ f, g \in V$.

From ii) in Lemma \ref{equivalentes} it now follows that $\|\cdot\|_V$ and $\|\cdot\|$ are
equivalent norms in $X(\mu)$.

Finally, consider $f, g \in X(\mu)$ such that $|f| \leq |g|$. By the monotonicity  of the norm
in $X(\mu)$, we obtain $ \|f\|_V = \big\|\,|f|\,\big\| \leq \big\|\,|g|\,\big\| = \|g\|_V$. \qed
\vspace{6pt}

\section{Some Consequences}\label{sec: conseq}

\subsection{Banach rectangular  f.s  with the subsequence property}

\begl  \label{cNRect} If $X(\mu)$ is a normed rectangular f.s., then $X(\mu)$ has an equivalent
strictly rectangular norm.
\terl
\noindent
Proof. Take $C >0 $ as in (\ref{dRectangular}). We define
\begeq \label{dNormaRect}
\|f\|_r \equS \sup \{ \|\chi_A f\|: A \in \Sig \}, \pt f \in X(\mu).
\tereq
Clearly $\|\cdot\|_r$ is a seminorm. Since
\begeq
\|f\|  \leq \|f\|_r \leq C \|f\|, \ \pt f \in X(\mu),
\tereq
it follows that $\|\cdot\|_r$ is an equivalent norm for $X(\mu)$.

Consider now $f \in X(\mu)$ and $A \in \Sig$. Then
\[
\|\chi_B \chi_A f \| = \| \chi_{B\inn A} f \| \leq \|f\|_r, \ \pt B \in \Sig.
\]
Hence $\|\chi_A f\|_r \leq \|f\|_r$. This shows that the norm $\|\cdot\|_r$ is strictly
rectangular.
\qed

\begc \label{Thm normaEFBEequiv} Let $X(\mu)$ be a Banach space of measurable functions.
Then, there is an equivalent norm under which $X(\mu)$  is  a Banach f.s. if, and only if,
$X(\mu)$ is a Banach rectangular f.s. with the subsequence property.
\terc

\noindent
Proof. Let us first assume that there is an equivalent norm $\|\cdot\|_B$ in $X(\mu)$ such
that   $X_B \equS (X(\mu),\|\cdot\|_B))$ is a Banach f.s. Then  Example \ref{prop BfsRectang}
indicates that $X(\mu)$ is Banach rectangular f.s. and from Corollary
\ref{BfsSubseqPpty} it follows that $X(\mu)$ has the subsequence property.

The remaining implication follows from Lemma  \ref{cNRect} and Theorem \ref{equivEFB}.
\qed \vspace{6pt}

As we will now see, applying Theorem \ref{equivEFB} we can obtain  for Banach rectangular
f.s. with the subsequence property, well known properties of Banach f.s.

\begd
Given a Banach rectangular f.s. with the subsequence property $X(\mu)$, we will denote by
$X(\mu)_V$  the Banach f.s. that we obtain when we consider in $X(\mu)$ the norm $\|\cdot
\|_V$ defined in (\ref{normaV}).
\terd

In the following we assume  $X(\mu)$ and $Y(\mu)$ to be  Banach rectangular f.s. with the
subsequence property. Naturally, a linear operator $T: X(\mu) \tien Y(\mu)$ is
\emph{positive} if $Tf \geq 0, \pt f \geq 0$.
By   Theorem \ref{equivEFB} the norms $\| \cdot\|$ are equivalent and $X(\mu)_V$ is a
Banach f.s. So $X(\mu)$ will have all those properties of $X(\mu)_V$ which only depend of
its topology. Next we give a simple example.

\begc
i) If $T: X(\mu) \tien Y(\mu)$ is a positive operator, then $T$ is continuous.
\espV
ii) If $X(\mu)$ is also a Banach rectangular space under the norm $\|\cdot\|_0$, then
$\|\cdot\|$ and $\|\cdot\|_0$ are equivalent norms.
\terc
\noindent
Proof. i) Clearly $T: X(\mu)_V \tien Y(\mu)_V$ is a positive operator. Being $X(\mu)_V$ and
$Y(\mu)_V$  B.f.s., it follows that $T$ is  bounded. Since the norm in $X(\mu)_V$
is equivalent with that of $X(\mu)$, and the norm in $Y(\mu)_V$ is equivalent with that of
$Y(\mu)$, this implies that $T: X(\mu) \tien Y(\mu)$ is bounded.

ii) Let $X(\mu)_j,$ $j=1,2,$ be the B.f.s. associated to the Banach rectangular f.s. $X(\mu)$ with the
norm $\|\cdot\|_j$. Then the map $I: X(\mu)_1 \tien X(\mu)_2$  given by $I(f) \equS f$ is a
linear isomorphism  and both $I$ and its inverse map are positive operators. It follows from
i) that both are bounded and this corresponds to the norms $\|\cdot\|_1$ and $\|\cdot\|_2$
being equivalent.\qed

\subsection{Normed sequence spaces}

\

We do not know if an  arbitrary Banach rectangular f.s. has the subsequence
property. However, we will  now show this holds for a rectangular sequence space.

We will say that  $X(\mu) = (\Ome,\Sig,\mu)$ is a \emph{normed  sequence space} (\emph{Banach
sequence space} when complete) when $\Ome:= \NN$,  $\Sig \equS 2^{\Ome}$ and   $\mu$ is  a
finite measure such that $\mu(\{n\}) > 0,  \pt n \in \NN$. The specific measure
will not be important. For definiteness we can take, for example, the measure defined by
$\mu (A) \equS \sum _{n \in A} \frac{1}{2^n}$.  Then any function $f: \NN \tien \RR$ is
measurable and, since the empty set is the only set in $\Sig$ with measure zero, each
equivalence class in $L^0(\mu)$ has only one member. So in this context we will identify
$L^0(\mu)$ with $L^0(\Sig)$ and also this last space with the set of real sequences
$\calS$. Let $A \cont \NN$, then considering the above identification, notice that
$\chi_A \equS \{a_n\}$, where $a_n = 1$ when $n \in A$ and $a_n = 0$ when $n \nin A$.

Let $X(\mu)$ be a normed sequence space. For each $n \in \NN$, we define
 the linear function $\pi_n: X(\mu) \tien \RR$ by $\pi_n(f) \equS f(n)$.

\begp\label{propRectSeq} If  $X(\mu)$ is a rectangular sequence space, then each linear
functional $\pi_n$ is continuous.  Hence $X(\mu)$ has the subsequence property.
\terp

\noindent
Proof. Take $C > 0$ as in 
Definition \ref{dRectangular}. Fix $n \in \NN$. If $g(n) = 0, \ \pt g \in X(\mu)$, the
conclusion is clear. So now assume that $g(n)\neq 0$, for some $g \in X(\mu)$. Then $g(n)
e_n = g \chi_{\{n\}} \in X(\mu)$. Hence $e_n \in X(\mu)$. Since $\|f(n) e_n\| = \| f e_n\|
\leq C \|f\|$, we have
\[
|\pi_n (f)| = |f(n)| \leq \frac{C \|f\|}{\|e_n\|}, \ \pt f \in X(\mu).
\]

It follows from the above inequality that $X(\mu)$ has the subsequence property. \qed

Using the above result together with  Theorem \ref{equivEFB}, we get the  following
characterization:

\begc \label{seqBFS}
Let  $X(\mu)\subset \calS $ be a  Banach sequence space.  Then $X(\mu)$ has an equivalent
norm under which it is a Banach  f.s. if and only if  $X(\mu)$ is a Banach rectangular f.s.
\terc

Let $X(\mu)$ be a Banach rectangular sequence space and, for simplicity,  assume that each
canonical sequence $e_n \in X(\mu)$. Naturally,  $e_n(k) \equS 1$ if $k = n$  and $e_n(k) =
0$ otherwise. Take $C > 0$ as in Definition \ref{dRectangular} and consider $1 \leq N < M$
and $A \equS \{1, \ldots, N\}$. Let $a_1, \ldots, a_{M} \in \RR$. Taking $s = (a_1, \ldots,
a_{M},  \ldots)$ we have
\[
\left\| \sum  _{n=1}^{N} a_n e_n\right\| = \|\chi_A s \| \leq C \|s\| =  C \left\| \sum
_{n=1}^{M} a_n e_n \right\|.
\]
By a well known result \cite[p.~411]{DJT}, this implies that $\{e_n\}$ is a basic sequence,
that is, it is a Schauder basis for the closure of its span, say $X_0$.

\begp $\{e_n\}$ is an unconditional basis for $X_0$.
\terp
\proof Let $s \equS \{a_n\} \in S(X_0)$, so $ x \equS \sum \subnui a_n e_n$ converges
in $X$. Take  $A \cont \NN$.  By Theorem \ref{equivEFB} multiplication by
$\chi_A$ is a bounded operator on $X(\mu)$.  Hence
$\chi_A x  = \sum \subnui \chi_A(n) a_n e_n$ and so $\chi_A s \in S(X_0)$.
Since  $X(\mu)$ has a rectangular norm,  from iii) of Theorem \ref{baseIncond}
we conclude that $\{e_n\}$ is an unconditional basis for $X_0$.
\qed

\end{document}